\documentclass{amsart}
\usepackage{amsfonts}
\usepackage{amssymb}
\usepackage{amsmath, amssymb, euscript,  enumerate}

\long\def\symbolfootnote[#1]#2{\begingroup\def\thefootnote{\fnsymbol{footnote}}\footnote[#1]{#2}\endgroup}

\newtheorem{theorem}{Theorem}

\newtheorem{corollary}{Corollary}

\begin{document}
\author{Ovidiu Munteanu}
\title{On the gradient estimate of Cheng and Yau}
\date{}

\begin{abstract}
We improve the well known local gradient estimate of Cheng and Yau in the
case when Ricci curvature has a negative lower bound.
\end{abstract}

\maketitle

\section{Introduction}

\symbolfootnote[0]{Reasearch partially supported by NSF grant No. DMS-1005484}

Let $M$ be an $m$ dimensional complete non-compact Riemannian manifold. The
local Cheng-Yau gradient estimate is a standard result in Riemannian
geometry, see \cite{CY}, also cf. \cite{LW}. It asserts that for $%
f:B_{p}(2R)\rightarrow \mathbb{R}$ harmonic and positive if the Ricci
curvature on $B_{p}(2R)$ has a lower bound $Ric\geq -(m-1)K$ for some $K\geq
0,$ then 
\begin{equation}
\sup_{B_{p}\left( R\right) }\left\vert \nabla \log f\right\vert \leq \left(
m-1\right) \sqrt{K}+\frac{C}{R}.  \label{a}
\end{equation}%
Notice that if $K=0$ it follows that a harmonic function with sublinear
growth on a manifold with non-negative Ricci curvature is constant. This
result is clearly sharp since on $\mathbb{R}^{n}$ there exist harmonic
functions which are linear.

When $K>0$, we can rescale the metric so that we take $K=1.$ In particular,
if $f$ is positive harmonic on $M$ and $Ric\geq -\left( m-1\right) $ on $M$
then it follows that $\left\vert \nabla \log f\right\vert \leq m-1$ on $M.$
This result is also sharp, in fact the equality case was recently
characterized in \cite{LW}.

This means that for $K=1$ in (\ref{a}) the factor $m-1$ on the right hand
side is sharp. However, the correction term that depends on the radius is
not sharp anymore. The purpose of this note is to prove a sharp version of
the local Cheng-Yau gradient estimate in the following form.

\begin{theorem}
Let $M$ be complete noncompact Riemannian manifold of dimension $m$ with
Ricci curvature bounded from below on a geodesic ball $B_{p}\left( 2R\right) 
$ by 
\begin{equation*}
Ric\geq -\left( m-1\right) .
\end{equation*}

There exists constants $C_{1}$ depending only on $m$ and $C_{2}$ universal
constant such that if $f:B_{p}\left( 2R\right) \rightarrow \mathbb{R}$ is a
positive harmonic function then 
\begin{equation*}
\sup_{x\in B_{p}\left( R\right) }\left\vert \nabla \log f\right\vert (x)\leq
m-1+\frac{C_{1}}{R}\exp \left( -C_{2}R\right) .
\end{equation*}
\end{theorem}

Certainly, our estimate becomes meaningful for large $R$. Let us also remark
that by rescaling the metric we obtain that if $Ric\geq -\left( m-1\right) K$
on $B_{p}\left( 2R\right) $ and $f$ is positive harmonic on $B_{p}\left(
2R\right) $ then 
\begin{equation*}
\sup_{B_{p}\left( R\right) }\left\vert \nabla \log f\right\vert \leq \left(
m-1\right) \sqrt{K}+\frac{C_{1}}{R}\exp \left( -C_{2}\sqrt{K}R\right) .
\end{equation*}%
The Cheng-Yau gradient estimate is a fundamental tool in geometric 
analysis. The technique has been applied in various situations, for an overview of the subject 
see e.g. \cite{SY, L}. In many applications the value of the constants in the estimate is quite important, 
for example the sharp version of the global estimate in the case when $Ric \ge -(m-1) $ on $M$ has been very instrumental
in rigidity theorems, see \cite{LW} and \cite{W}. 
The importance of our Theorem, besides being sharp, is that the estimate does not diverge when integrated on minimizing 
geodesics. For example, in Corollary 1 in this paper we have established a sharp lower bound for the Green's
function on a nonparabolic manifold with a negative lower bound for Ricci curvature. Similar sharp upper or lower 
bound estimates can be established for harmonic functions defined on manifolds with boundary.

\section{Proof of the result}

\begin{proof}
The proof of the Theorem follows the standard argument of Cheng and Yau. We
apply the Bochner technique to $\phi ^{2}\left\vert \nabla \log f\right\vert
^{2},$ where $\phi $ is a cut-off function with support in $B_{p}\left(
R\right) .$ The difference here is that we use a judicious choice of cut-off 
$\phi $ which improves the argument.

We will prove the following statement: Let $f:B_{p}\left( 2R\right)
\rightarrow \mathbb{R}$ be positive and harmonic, where $R\geq R_{0}\left(
m\right) $ with $R_{0}$ depending only on $m$. 

Assume that $Ric\geq -\left(
m-1\right) $ on $B_{p}\left( 2R\right) ,$ then 
\begin{equation}
\left\vert \nabla \log f\right\vert \left( p\right) \leq m-1+C_{1}\exp
\left( -C_{2}R\right) .  \label{*}
\end{equation}%
It is evident that (\ref{*}) proves Theorem 1, because we can apply (\ref{*}%
) for each $x\in B_{p}\left( R\right) $. Of course, since (\ref{*}) is true
for any radius, we can apply it for the radius $\frac{R}{2}$, with $R$ the
same as in Theorem 1.

Let $h=\log f$ . Consider also the function $\phi :[0,\infty )\rightarrow 
\mathbb{R}$ defined by 
\begin{eqnarray*}
\phi \left( 0\right) &=&1, \\
\phi \left( r\right) &=&1-\exp \alpha \left( R-\frac{R^{2}}{r}\right) \ \ \ 
\text{for\ \ }0<r\leq R \\
\phi \left( r\right) &=&0\ \ \ \text{for\ \ }r>R.
\end{eqnarray*}%
where $\alpha >0$ is a (small) number that will be set later.

We compute directly that for $r<R,$ 
\begin{eqnarray}
\phi ^{\prime }\left( r\right) &=&-\alpha \frac{R^{2}}{r^{2}}\exp \alpha
\left( R-\frac{R^{2}}{r}\right) ,  \label{b} \\
\phi ^{\prime \prime }\left( r\right) &=&\left( 2\alpha \frac{R^{2}}{r^{3}}%
-\alpha ^{2}\frac{R^{4}}{r^{4}}\right) \exp \alpha \left( R-\frac{R^{2}}{r}%
\right) ,  \notag \\
\phi ^{\prime }\left( 0\right) &=&\phi ^{\prime \prime }\left( 0\right) =0. 
\notag
\end{eqnarray}%
By applying $\phi $ to the distance function from $p$ we obtain a cut-off
function on $M,$ with support in $B_{p}\left( R\right) $, which we continue
to denote with $\phi ,$ i.e. $\phi \left( x\right) =\phi \left( r\left(
x\right) \right) $.

Let now consider $G:B_{p}\left( R\right) \rightarrow \mathbb{R},$ 
\begin{equation*}
G=\phi ^{2}\left\vert \nabla h\right\vert ^{2}.
\end{equation*}%
Since $G$ is nonnegative on $B_{p}\left( R\right) $ and $G=0$ on $\partial
B_{p}\left( R\right) ,$ it follows that $G$ attains a maximum point in the
interior of $B_{p}\left( R\right) .$ Let $x_{0}$ be this maximum point. We
will assume from now on that $r\left( x\right) =d\left( p,x\right) $ is
smooth at $x_{0},$ so $\phi $ and $G$ are also smooth at $x_{0}.$ If this is
not the case, then one can use a support function at $x_{0}$, which is
smooth. The computations will be similar and will imply the same result. In
this standard argument, it is essential that $\phi \left( r\right) $ be
nonincreasing, which is true in our case, see \cite{C}.

Hence at $x_{0}$ we have, by the maximum principle: 
\begin{eqnarray}
\left\vert \nabla G\left( x_{0}\right) \right\vert  &=&0  \label{b'} \\
\Delta G\left( x_{0}\right)  &\leq &0.  \notag
\end{eqnarray}%
The following part of the argument is well known, but we include it for
completeness. By the Bochner formula, since $h=\log f$ and $\Delta
h=-\left\vert \nabla h\right\vert ^{2},$ we get

\begin{eqnarray}
\frac{1}{2}\Delta \left\vert \nabla h\right\vert ^{2} &=&\left\vert
h_{ij}\right\vert ^{2}+\left\langle \nabla h,\nabla \left( \Delta h\right)
\right\rangle +Ric\left( \nabla h,\nabla h\right)  \notag \\
&\geq &\left\vert h_{ij}\right\vert ^{2}-\left\langle \nabla h,\nabla
\left\vert \nabla h\right\vert ^{2}\right\rangle -\left( m-1\right)
\left\vert \nabla h\right\vert ^{2}.  \label{c}
\end{eqnarray}%
Moreover, choosing an ortonormal frame $\left\{ e_{i}\right\} _{i}$ such
that $e_{1}=\frac{\nabla h}{\left\vert \nabla h\right\vert },$ we have

\begin{gather*}
\left\vert h_{ij}\right\vert ^{2}\geq \left\vert h_{11}\right\vert
^{2}+\sum\limits_{\alpha >1}\left\vert h_{\alpha \alpha }\right\vert
^{2}+2\sum\limits_{\alpha >1}\left\vert h_{1\alpha }\right\vert ^{2} \\
\geq \left\vert h_{11}\right\vert ^{2}+2\sum\limits_{\alpha >1}\left\vert
h_{1\alpha }\right\vert ^{2}+\frac{1}{m-1}\left\vert \sum_{\alpha
>1}h_{\alpha \alpha }\right\vert ^{2} \\
=\left\vert h_{11}\right\vert ^{2}+2\sum\limits_{\alpha >1}\left\vert
h_{1\alpha }\right\vert ^{2}+\frac{1}{m-1}\left\vert \left\vert \nabla
h\right\vert ^{2}+h_{11}\right\vert ^{2} \\
\geq \frac{m}{m-1}\left( \left\vert h_{11}\right\vert
^{2}+\sum\limits_{\alpha >1}\left\vert h_{1\alpha }\right\vert ^{2}\right) +%
\frac{1}{m-1}\left\vert \nabla h\right\vert ^{4}+\frac{2}{m-1}%
h_{11}\left\vert \nabla h\right\vert ^{2}.
\end{gather*}

On the other hand, notice that 
\begin{eqnarray*}
\left\langle \nabla \left\vert \nabla h\right\vert ^{2},\nabla
h\right\rangle  &=&2h_{ij}h_{i}h_{j}=2h_{11}\left\vert \nabla h\right\vert
^{2}\ \ \text{and} \\
\left\vert \nabla \left\vert \nabla h\right\vert ^{2}\right\vert ^{2}
&=&4\left\vert h_{ij}h_{j}\right\vert ^{2}=4h_{1i}^{2}\left\vert \nabla
h\right\vert ^{2}
\end{eqnarray*}%
which imply 
\begin{equation*}
\left\vert h_{ij}\right\vert ^{2}\geq \frac{m}{4\left( m-1\right) }%
\left\vert \nabla h\right\vert ^{-2}\left\vert \nabla \left\vert \nabla
h\right\vert ^{2}\right\vert ^{2}+\frac{1}{m-1}\left\vert \nabla
h\right\vert ^{4}+\frac{1}{m-1}\left\langle \nabla \left\vert \nabla
h\right\vert ^{2},\nabla h\right\rangle .
\end{equation*}

Then using this in (\ref{c}) we conclude that 
\begin{gather}
\frac{1}{2}\Delta \left\vert \nabla h\right\vert ^{2}\geq \frac{m}{4\left(
m-1\right) }\left\vert \nabla h\right\vert ^{-2}\left\vert \nabla \left\vert
\nabla h\right\vert ^{2}\right\vert ^{2}+\frac{1}{m-1}\left\vert \nabla
h\right\vert ^{4}  \label{d} \\
-\frac{m-2}{m-1}\left\langle \nabla \left\vert \nabla h\right\vert
^{2},\nabla h\right\rangle -\left( m-1\right) \left\vert \nabla h\right\vert
^{2}  \notag
\end{gather}

By (\ref{d}) it results%
\begin{gather*}
\frac{1}{2}\Delta G\geq \frac{m}{4\left( m-1\right) }\phi
^{4}G^{-1}\left\vert \nabla \left( \phi ^{-2}G\right) \right\vert ^{2}+\frac{%
1}{m-1}\phi ^{-2}G^{2} \\
-\frac{m-2}{m-1}\phi ^{2}\left\langle \nabla \left( \phi ^{-2}G\right)
,\nabla h\right\rangle -\left( m-1\right) G+\frac{1}{2}\phi ^{-2}\left(
\Delta \phi ^{2}\right) G+\left\langle \nabla \phi ^{2},\nabla \left( \phi
^{-2}G\right) \right\rangle .
\end{gather*}%
At the maximum point $x_{0}$ it follows from (\ref{b'}) that 
\begin{eqnarray*}
0 &\geq &\frac{1}{m-1}G^{2}-\left( m-1\right) \phi ^{2}G+\frac{m}{m-1}%
\left\vert \nabla \phi \right\vert ^{2}G+\frac{2\left( m-2\right) }{m-1}\phi
\left\langle \nabla \phi ,\nabla h\right\rangle G \\
&&+\frac{1}{2}\left( \Delta \phi ^{2}\right) G-4\left\vert \nabla \phi
\right\vert ^{2}G.
\end{eqnarray*}

Since $\phi \left\langle \nabla \phi ,\nabla h\right\rangle \geq -\left\vert
\nabla \phi \right\vert G^{\frac{1}{2}}$ we infer from the above inequality
that:%
\begin{eqnarray*}
0 &\geq &\frac{1}{2}\left( m-1\right) \left( \Delta \phi ^{2}\right)
G-\left( 3m-4\right) \left\vert \nabla \phi \right\vert ^{2}G-\left(
m-1\right) ^{2}\phi ^{2}G \\
&&-2\left( m-2\right) \left\vert \nabla \phi \right\vert G^{\frac{3}{2}%
}+G^{2}.
\end{eqnarray*}%
This can be written as 
\begin{equation}
-\left( m-1\right) \phi \Delta \phi +\left( 2m-3\right) \left\vert \nabla
\phi \right\vert ^{2}+\left( m-1\right) ^{2}\phi ^{2}+2\left( m-2\right)
\left\vert \nabla \phi \right\vert G^{\frac{1}{2}}\geq G.  \label{e}
\end{equation}%
We point out that (\ref{e}) is true for any cut-off $\phi ,$ but now we want
to estimate the left hand side of (\ref{e}) from above using our choice of
cut-off $\phi $. In existing arguments in the literature one would bound $%
\left( m-1\right) ^{2}\phi ^{2}\leq \left( m-1\right) ^{2}$ and then deal
with the terms involving $\left\vert \nabla \phi \right\vert $ and $\Delta
\phi .$ Our strategy is to make use of the fact that $\phi $ is small near
the boundary. So, if $x_{0}$ is close to $\partial B_{p}\left( R\right) $
then it will be easy to see that $G\left( x_{0}\right) \leq \left(
m-1\right) ^{2}.$ On the other hand, if that is not true, then $\left\vert
\nabla \phi \right\vert $ and $\Delta \phi $ will be very small for our
choice of $\phi ,$ because $\phi $ decays very slowly until the boundary. 

First, observe that using standard local gradient estimate (\ref{a}) we can
estimate 
\begin{equation*}
\sup_{B_{p}\left( R\right) }\left\vert \nabla h\right\vert \leq \left(
m-1\right) +\frac{c_{1}}{R},
\end{equation*}%
where $c_{1}$ is a constant depending only on $m.$ Therefore, choosing $%
R_{0}\left( m\right) $ sufficiently big, we can guarantee that $\left(
m-1\right) +\frac{c_{1}}{R}\leq m.$ Consequently, 
\begin{equation*}
G^{\frac{1}{2}}\left( x_{0}\right) =\left\vert \nabla h\right\vert \left(
x_{0}\right) \phi \left( x_{0}\right) \leq \sup_{B_{p}\left( R\right)
}\left\vert \nabla h\right\vert \leq m.
\end{equation*}%
Next, let us consider the function $\theta :[0,R]\rightarrow \mathbb{R},$ 
\begin{equation*}
\theta \left( r\right) =\frac{R^{2}}{r^{2}}\exp \alpha \left( R-\frac{R^{2}}{%
r}\right) .
\end{equation*}%
Then we see that 
\begin{equation*}
\theta ^{\prime }\left( r\right) =\left( \alpha \frac{R^{4}}{r^{4}}-2\frac{%
R^{2}}{r^{3}}\right) \exp \alpha \left( R-\frac{R^{2}}{r}\right) ,
\end{equation*}%
so by assuming $R_{0}>\frac{2}{\alpha }$ it follows that $\theta $ is
increasing on $[0,R].$ Consequently, 
\begin{equation*}
\theta \left( r\right) \leq \theta \left( R\right) =1.
\end{equation*}%
By (\ref{b}), we have thus proved that 
\begin{equation}
\left\vert \nabla \phi \right\vert \leq \alpha .  \label{f}
\end{equation}%
Furthermore, the Laplacian comparison theorem states that 
\begin{equation*}
\Delta r\left( x_{0}\right) \leq \left( m-1\right) \left( 1+\frac{1}{r\left(
x_{0}\right) }\right) ,
\end{equation*}%
so that we find 
\begin{eqnarray*}
-\Delta \phi  &=&-\phi ^{\prime }\Delta r-\phi ^{\prime \prime } \\
&\leq &\left( \alpha \left( m-1\right) \frac{R^{2}}{r^{2}}\left( 1+\frac{1}{r%
}\right) +\left( -2\alpha \frac{R^{2}}{r^{3}}+\alpha ^{2}\frac{R^{4}}{r^{4}}%
\right) \right) \exp \alpha \left( R-\frac{R^{2}}{r}\right)  \\
&=&\alpha \left( \left( m-1\right) \frac{R^{2}}{r^{2}}+\left( m-3\right) 
\frac{R^{2}}{r^{3}}+\alpha \frac{R^{4}}{r^{4}}\right) \exp \alpha \left( R-%
\frac{R^{2}}{r}\right) .
\end{eqnarray*}%
Consequently, since $0\leq \phi \leq 1$ we get%
\begin{equation*}
-\phi \Delta \phi \leq \alpha \left( \left( m-1\right) \frac{R^{2}}{r^{2}}%
+\left( m-3\right) \frac{R^{2}}{r^{3}}+\alpha \frac{R^{4}}{r^{4}}\right)
\exp \alpha \left( R-\frac{R^{2}}{r}\right) .
\end{equation*}

Plugging all these into (\ref{e}) and using (\ref{b}) and (\ref{f}) it
follows that at $x_{0}$ 
\begin{eqnarray*}
G &\leq &\left( m-1\right) ^{2}\phi -\left( m-1\right) \phi \Delta \phi
+\left( 2m-3\right) \alpha \left\vert \nabla \phi \right\vert +2m\left(
m-2\right) \left\vert \nabla \phi \right\vert  \\
&\leq &\left( m-1\right) ^{2}\left( 1-\exp \alpha \left( R-\frac{R^{2}}{r}%
\right) \right)  \\
&&+\{\alpha \left( 2m\left( m-2\right) +\left( 2m-3\right) \alpha +\left(
m-1\right) ^{2}\right) \frac{R^{2}}{r^{2}} \\
&&+\alpha \left( m-1\right) \left( m-3\right) \frac{R^{2}}{r^{3}}+\left(
m-1\right) \alpha ^{2}\frac{R^{4}}{r^{4}}\}\exp \alpha \left( R-\frac{R^{2}}{%
r}\right) .
\end{eqnarray*}%
We can simplify this using that 
\begin{equation*}
\frac{R^{2}}{r^{2}}\leq \frac{R^{4}}{r^{4}},\ \ \frac{R^{2}}{r^{3}}\leq 
\frac{R^{4}}{r^{4}},
\end{equation*}%
and furthermore that for  $\alpha \leq 2/3,$%
\begin{equation*}
2m\left( m-2\right) +\left( 3m-4\right) \alpha +\left( m-1\right)
^{2}+\left( m-1\right) (m-3)\leq 4\left( m-1\right) ^{2}.
\end{equation*}

It results that at $x_{0}$ we have: 
\begin{equation}
G\leq \left( m-1\right) ^{2}+\left( m-1\right) ^{2}\left( 4\alpha \frac{R^{4}%
}{r^{4}}-1\right) \exp \alpha \left( R-\frac{R^{2}}{r}\right) .  \label{g}
\end{equation}%
Let us denote now $u:[0,R]\rightarrow \mathbb{R},$%
\begin{equation*}
u\left( r\right) =\left( 4\alpha \frac{R^{4}}{r^{4}}-1\right) \exp \alpha
\left( R-\frac{R^{2}}{r}\right) .
\end{equation*}%
We want to find the maximum value of $u$ on $[0,R],$ which by (\ref{g}) will
imply the desired estimate for $\left\vert \nabla h\right\vert \left(
p\right) .$

We have the following two cases:

Case 1: $r\geq \left( 4\alpha \right) ^{\frac{1}{4}}R$

Then it is easy to see that $u\left( r\right) \leq 0.$

Case 2: $r<\left( 4\alpha \right) ^{\frac{1}{4}}R$

In this case it follows trivially that 
\begin{equation*}
u\left( r\right) \leq 4\alpha \frac{R^{4}}{r^{4}}\exp \alpha \left( R-\frac{%
R^{2}}{r}\right) .
\end{equation*}%
It is then easy to see that the function $w:[0,R]\rightarrow \mathbb{R}$%
\begin{equation*}
w\left( r\right) =\frac{R^{4}}{r^{4}}\exp \alpha \left( R-\frac{R^{2}}{r}%
\right)
\end{equation*}%
is increasing on $[0,R]$ if we assumed $R_{0}\geq \frac{4}{\alpha }.$

Indeed, 
\begin{equation*}
w^{\prime }(r)=\left( \alpha \frac{R^{6}}{r^{6}}-4\frac{R^{4}}{r^{5}}\right)
\exp \alpha \left( R-\frac{R^{2}}{r}\right) \geq 0
\end{equation*}%
if 
\begin{equation*}
r\leq \frac{\alpha }{4}R^{2}.
\end{equation*}%
But this is clearly true using that $R\geq R_{0}\geq \frac{4}{\alpha }$ and
that $r\leq R.$

We have proved that $w$ is increasing on $\left[ 0,R\right] ,$ therefore for
any $r<\left( 4\alpha \right) ^{1/4}R$ we get%
\begin{eqnarray*}
u\left( r\right)  &\leq &4\alpha w\left( r\right) \leq 4\alpha w\left(
\left( 4\alpha \right) ^{1/4}R\right)  \\
&=&\exp \left( -\alpha \left( \left( 4\alpha \right) ^{-1/4}-1\right)
R\right)  \\
&=&\exp \left( -C_{2}R\right) .
\end{eqnarray*}%
Note that $\alpha $ can be chosen a number not depending on $m$, and then $%
C_{2}$ is independent on $m,$ too. For example, we can take $\alpha =2^{-6}$
and then $C_{2}=2^{-6}.$

Summing up, based on Case 1 and Case 2, it results that for $R\geq {R_{0}}$ 
\begin{eqnarray*}
\left\vert \nabla \log f\right\vert ^{2}\left( p\right) &\leq &G\left(
x_{0}\right) \\
&\leq &\left( m-1\right) ^{2}+\left( m-1\right) ^{2}\exp \left(
-C_{2}R\right) ,
\end{eqnarray*}%
which is what we claimed.

This proves the local gradient estimate, when the distance function from $p$
is smooth at $x_{0},$ the maximum point of $G$. As mentioned in the
beginning, if this is not the case, we consider $\gamma \left( t\right) $,
the minimizing geodesic from $p$ to $x_{0}.$ Certainly, we may assume that $%
x_{0}\neq p,$ since otherwise the estimates are trivial based on the
definition of $\phi .$

Then take $q=\gamma \left( \varepsilon \right) $ for $\varepsilon $ small
and define 
\begin{equation*}
\psi \left( x\right) :=\phi \left( d\left( q,x\right) +\varepsilon \right) .
\end{equation*}

Since $x$ is in the cut locus of $p,$ then it is not in the cut locus of $q,$
so $\psi $ is smooth at $x_{0}.$ Moreover, we have $d\left( q,x\right)
+\varepsilon \geq r\left( x\right) $ for any $x\in M$ and $d\left(
q,x_{0}\right) +\varepsilon =r\left( x_{0}\right) .$ Therefore, using that $%
\phi $ is decreasing on $\mathbb{R}$ it follows that $\psi \left( x\right)
\leq \phi \left( x\right) $ for any $x\in M$ and $\psi \left( x_{0}\right)
=\phi \left( x_{0}\right) .$

This means that $x_{0}$ is the maximum point of $\widetilde{G}:=\psi
^{2}\left\vert \nabla h\right\vert ^{2}$, which now is smooth at $x_{0}.$
Performing all the above computations and letting $\varepsilon \rightarrow 0$
it is not difficult to see that we still obtain (\ref{g}) at $x_{0}.$ The
rest of the argument is the same.

This proves (\ref{*}), which as explained in the beginning of the proof also
proves the Theorem. 
\end{proof}

As a consequence, we get the following sharp lower bound for the Green's
function, by integrating the estimate in Theorem 1 along minimizing
geodesics. We call $M$ nonparabolic if it admits a positive symmetric
Green's function.

\begin{corollary}
Let $M$ be complete noncompact Riemannian manifold of dimension $m$ with $%
Ric\geq -\left( m-1\right) $ on $M.$ If $M$ is nonparabolic, then there
exists a constant $C$ depending only on $m$ such that 
\begin{equation*}
\sup_{x\in \partial B_{p}\left( R\right) }G\left( p,x\right) \geq
C\inf_{y\in \partial B_{p}\left( 1\right) }G\left( p,y\right) \cdot
e^{-\left( m-1\right) R},
\end{equation*}%
where $G\left( p,x\right) $ is the positive symmetric Green's function with
a pole at $p\in M$.
\end{corollary}

On $\mathbb{H}^{n}$ let us take now the Green's function $f\left( x\right)
:=G\left( q,x\right) $ with a pole at $q.$ It is known that we have the
following formula%
\begin{equation}
f\left( x\right) =\int_{d\left( q,x\right) }^{\infty }\frac{dt}{A\left(
t\right) }  \label{h}
\end{equation}%
where $A\left( t\right) $ is the area of $\partial B_{q}\left( t\right) .$
Then it can be checked that for $p$ such that $d\left( p,q\right) =R\geq 1$ 
\begin{equation*}
\left( n-1\right) +C^{-1}e^{-2R}\leq \left\vert \nabla \log f\right\vert
\left( p\right) \leq \left( n-1\right) +Ce^{-2R}.
\end{equation*}%
This shows that the exponential-type decay in our Theorem is sharp. Also,  (%
\ref{h}) shows that the decay estimate in the Corollary is sharp.

{\small DEPARTMENT OF MATHEMATICS, COLUMBIA UNIVERSITY }

{\small NEW YORK, NY 10027}\newline
{\small E-mail address: omuntean@math.columbia.edu}

\end{document}